\documentclass[12pt]{article}
\usepackage{amssymb}
\usepackage[dvipdfm]{graphicx}
\usepackage{amsfonts}
\usepackage{latexsym}
\usepackage{amsthm}
\usepackage{amsmath}

\usepackage{amssymb, amscd}

\usepackage{url}

\newtheorem{problem}{Problem}
\newtheorem{theo}[problem]{Theorem}
\newtheorem{rem}[problem]{Remark}
\newtheorem{prob}[problem]{Problem}
\newtheorem{defin}[problem]{Definition}
\newtheorem{prop}[problem]{Proposition}
\newtheorem{cor}[problem]{Corollary}

\newtheorem{exam}[problem]{Example}

\begin{document}

 \title{{Chessboard complexes indomitable}\footnote{This paper is an expanded version of our preprint \cite{X}, with added
 Theorem~\ref{thm:general} and its consequences.}}

\author{Sini\v sa T.\ Vre\' cica\\ {\small Faculty of Mathematics}\\[-2mm] {\small University of Belgrade}
\\[-2mm] {\small vrecica$@$matf.bg.ac.rs}
 \and Rade  T.\ \v Zivaljevi\' c\\ {\small Mathematical Institute}\\[-2mm] {\small SASA, Belgrade}\\[-2mm]
 {\small rade$@$mi.sanu.ac.rs} }
\date{April 10, 2011}

\maketitle  \vspace{-1cm}

\begin{abstract}
We give a simpler, degree-theoretic proof of the striking new
Tverberg type theorem of Blagojevi\' c, Ziegler and Matschke,
arXiv:0910.4987v2. Our method also yields some new examples of
``constrained Tverberg theorems'' including a simple colored
Radon's theorem for $d+3$ points in $\mathbb{R}^d$. This gives us
an opportunity to review some of the highlights of this beautiful
theory and reexamine the role of chessboard complexes in these and
related problems of topological combinatorics.
\end{abstract}
\renewcommand{\thefootnote}{$\ast$} \footnotetext{Supported by Grants 144014
and 144026 of the Serbian Ministry of Science and Technology.}

\section{Introduction}
Chessboard (simplicial) complexes and their relatives have been
for decades an important theme of topological combinatorics with
often unexpected applications in group theory, representation
theory, commutative algebra, Lie theory, computational geometry,
algebraic topology, and combinatorics, see \cite{Ata04},
\cite{AuFie07}, \cite{BLVZ}, \cite{FrHa98}, \cite{Ga79},
\cite{Jo07}, \cite{ShaWa04}, \cite{VZ94}, \cite{VZ07},
\cite{Wa03}, \cite{Zie94}, \cite{ZV92}. The books \cite{Jo-book}
and \cite{M}, as well as the review papers \cite{Wa03} and
\cite{Z04}, cover selected topics of the theory of chessboard
complexes and contain a more complete list of related
publications.

\medskip
Chessboard complexes originally appeared in \cite{Ga79} as coset
complexes of the symmetric group, closely related to Coxeter  and
Tits coset complexes. In combinatorics they appeared as
``complexes of partial injective functions'' \cite{ZV92},
``multiple deleted joins'' \cite{ZV92}, complexes of all partial
matchings in complete bipartite graphs, and the complexes of all
non-taking rook configurations \cite{BLVZ}.

Recently a naturally defined subcomplex of the chessboard complex,
referred to as the ``cycle-free chessboard complex'', has emerged
in the context of stable homotopy theory (\cite{AuFie07} and
\cite{Fie07}), where it was introduced as a tool for evaluating
the symmetric group analogue for the cyclic homology of algebras.

\medskip
In our own research \cite{ZV92, VZ94} chessboard complexes
appeared as a tool for the resolution of the well known {\em
colored Tverberg problem}, see \cite{M} and \cite{Z04} for the
history of the problem and its connections with other well known
problems of discrete and computational geometry. In these papers
the fundamental role of chessboard complexes for colored Tverberg
type problems was discovered, and the importance of Borsuk-Ulam
type questions for equivariant maps defined on joins of chessboard
complexes recognized.

\medskip
Next fifteen years witnessed little progress and it is probably
safe to say that majority of specialists, including ourselves,
arrived at a conclusion that the limits of the method are reached
and a new progress towards better bounds in the colored Tverberg
problem difficult to expect.

\medskip
Consequently it was indeed a wonderful surprise when Pavle
Blagojevi\' c and G\" unter Ziegler \cite{B-Z} (and Benjamin
Matschke, see \cite{BMZ}, the second version of \cite{B-Z}) proved
the opposite and established so far the most natural and elegant
version of (type~A) colored Tverberg theorem.

\medskip
Motivated by the breakthrough of Blagojevi\' c and Ziegler we
prove some new Borsuk-Ulam type results for joins of chessboard
complexes (Propositions~\ref{prop:degree}, \ref{prop:posledica}
and Theorem~\ref{thm:general}) leading to a shorter and
conceptually simpler proof\footnote{After the preliminary version
\cite{X} of our paper was released  and shared with a circle of
specialists, we were kindly informed by P.~Blagojevi\' c that
B.~Matschke has also discovered a proof that simplifies their
original approach. This proof is incorporated in \cite{BMZ-2}.} of
their main result. Among other consequences of our approach are
new constrained Tverberg theorems close in spirit to results of
S.~Hell \cite{H}, the simplest example being the ``Colored Radon's
theorem'' (Corollary~\ref{cor:Radon}).

\section{Preliminaries}

\subsection{Chessboard complexes $\Delta_{m,n}$ }\label{sec:coloring}

A function $f : A\rightarrow B$ can be interpreted as a {\em
labeling} of elements of a set $A$ by labels from a set $B$. A
partial labeling of $A$ is a function $\phi : D\rightarrow B$
where $D=D(\phi)$ the domain of $\phi$ is a subset of $A$. By
convention a partial function $\phi$ from $A$ to $B$ is often
identified with its graph $\Gamma(\phi) = \{(i,\phi(i))\mid i\in
D(\phi)\}\subset A\times B$. It follows that if $A$ and $B$ are
finite then the set ${\rm L}(A,B)$ of all (partial) labelings of
$A$ by labels from $B$ is a simplicial complex on $A\times B$ as a
set of vertices. If the cardinalities of sets $A$ and $B$ are
respectively $m$ and $n$ then the complex ${\rm L}(A,B)$ is
immediately recognized as a complex isomorphic to the join
$A^{\ast n} \cong [m]^{\ast n}$ of $n$ copies of an $m$-element
set ($0$-dimensional complex) $A$.

\medskip
If we restrict our attention to {\em partial injective functions}
we immediately arrive at the chessboard complex
$$
\Delta_{m,n}:=\Delta_{A,B} = \{ \Gamma(\phi)\in {\rm L}(A,B)\mid
\phi \mbox{ {\rm is an injective function}}\}
$$
in the form it was introduced in \cite{ZV92, VZ94} in the context
of the colored Tverberg problem. The name ``chessboard complex''
is motivated by the fact that each simplex $\Gamma(\phi)$ can be
interpreted as a non-taking rook placement on a $(m\times
n)$-chessboard $A\times B\cong [m]\times [n]$, i.e.\ a
configuration of rooks where no rook threatens any other.

\medskip
The ``injectivity restriction'' on a labeling is quite natural
from the point of view of ``constrained labelings'' in the sense
of Hell \cite{H}, where only some of the labelings are admissible.
For example if $A$ is the vertex set of a graph $\Gamma = (A,E)$,
then it is natural to ask that $\phi$ is a proper labeling in the
sense that two adjacent vertices always receive different labels.

If $A=[m], B=[n]$ and $\Gamma = K_m$ is the complete graph, then
the simplicial complex of all admissible (partial) labelings of
$A$ by $n$ distinct labels is precisely the chessboard complex
$\Delta_{m,n}$. More generally if $\Gamma$ is a disjoint union of
cliques, i.e.\ if there is a partition $A=A_1\cup\ldots\cup A_k$
such that $\Gamma = K_{A_1}\cup\ldots\cup K_{A_k}$ is the union of
complete graphs, then the complex of admissible labelings is
isomorphic to the join
\begin{equation}\label{eqn:chess-join}
\Delta_{A_1,B}\ast\ldots\ast \Delta_{A_k,B}
\end{equation}
of chessboard complexes.

\subsection{Colored Tverberg problem}
\label{sec:CTP}

\medskip
Suppose that $C\subset \mathbb{R}^d$ is a finite set and let $\psi
: C \rightarrow [k+1]$ be a coloring of this set
(Section~\ref{sec:coloring}) by $k+1$ colors where $k\leq d$. The
coloring is always supposed to be a {\em strict coloring} in the
sense that the coloring function $\psi$ is an epimorphism, i.e.\
that all listed colors are used.

A subset $X\subset C$ is called {\em multicolored} if the
restriction of $\psi$ on $X$ is injective, i.e.\ if elements of
$X$ are all colored by different colors. In this case the
(possibly degenerate) simplex $\sigma := {\rm conv}(X)$ is also
referred to as a multicolored or a {\em rainbow} simplex.

\begin{enumerate}
\item[$\bullet$] Informally speaking the Colored Tverberg problem
asks for conditions on the coloring function $\psi$ which
guarantee the existence of many, vertex-disjoint rainbow simplices
which have a nonempty intersection.
\end{enumerate}
If the colors are evenly distributed we are led to the following
version of the general problem in the form it was recorded in
$\cite{Z04}$.
\begin{prob}\label{prob:color}
For given integers $r,k,d$ such that $k\leq d$, determine the
minimum number $t = t(r,k,d)$ such that for any $C\subset
\mathbb{R}^d$ of size $t(k+1)$ and each strict coloring $\psi :
C\rightarrow [k+1]$ such that $C_j:=\psi^{-1}(j)$ has $t$
elements, there exist $r$ disjoint, multicolored sets $B_j\subset
C$ such that
\begin{equation}\label{eqn:presek}
\bigcap_{j=1}^r {\rm conv}(B_j)\neq\emptyset.
\end{equation}
\end{prob}
The reader is again referred to the book \cite{M} and reviews
\cite{Z04, User2} for a much more detailed presentation with a
fairly complete set of references. Here we recall only that the
problem of evaluating $t(r,d,d)$ was originally proposed by B\'
ar\' any and Larman \cite{BL}, after it was observed by
B\'{a}r\'{a}ny, F\"{u}redy and Lov\'{a}sz \cite{BFL} that the
``weak colored Tverberg theorem'' $t(r,d,d)<+\infty$ resolves a
number of interesting conjectures in discrete and computational
geometry (halving hyperplanes problem, point selection problem,
weak $\epsilon$-net problem, hitting set problem).

After the preliminary result \cite{BL} that $t(r,2,2)=r$ and
$t(2,d,d)=2$, the general bounds $t(r,d,d)\leq 2r-1$ and
$t(r,k,d)\leq 2r-1$ were established in \cite{ZV92}, respectively
\cite{VZ94} for all primes $r$. Subsequently \cite{User2}, and
without introducing really new ideas, the result was extended to
the case of prime powers. Note that the distinction between the
cases $k=d$ and $k<d$ is important since in the latter case (for
dimensional reasons) there is an additional constraint $r\leq
d/(d-k)$. For this reason the case $k=d$ is referred to as the
``type A colored Tverberg problem'' while the case $k<d$ is known
as the ``type B colored Tverberg problem'', see \cite{Z04}.

Note that the bound $t(r,k,d)\leq 2r-1$ for $k<d$ was shown in
\cite{VZ94} to be tight if $r$ is a prime, while the central new
result of \cite{B-Z} (confirming the conjecture from \cite{BL}) is
the equality $t(r,d,d)=r$ provided $p=r+1$ is a prime number.

\subsection{Topology enters the scene}

The essence of the original breakthrough \cite{ZV92, VZ94},
leading to the inequality $t(r,k,d)\leq 2r-1$, was the observation
that the colored Tverberg problem is closely related to a question
of Borsuk-Ulam type for joins of chessboard complexes. More
precisely it was shown that both type A and type B cases of the
problem follow from the nonexistence of a
$\mathbb{Z}/r$-equivariant map
\begin{equation}\label{eqn:nasa}
(\Delta_{r,2r-1})^{\ast (k+1)}\longrightarrow S(W_r^{\oplus d})
\end{equation}
where $W_r$ is the standard $(r-1)$-dimensional real
representation of the cyclic group $\mathbb{Z}/r$, i.e.\ the
representation obtained by removing the trivial from the regular
representation of $\mathbb{Z}/r$.

Indeed, let $C\subset \mathbb{R}^d$ be a set of size $(2r-1)(k+1)$
and let $\psi : C\rightarrow [k+1]$ be a coloring function such
that each $C_j :=\psi^{-1}(j)$ has $2r-1$ elements. We are
supposed to show that there exist pairwise disjoint subsets
$B_j,\, j=1,\ldots, r,$ satisfying (\ref{eqn:presek}) such that
the restriction $\psi\vert_{B_j}$ is injective for each $j$. Let
$\phi$ be a partial labeling of $C$ (Section~\ref{sec:coloring})
such that $D(\phi)=\cup_{j=1}^r~B_j$ and $B_j := \phi^{-1}(j)$.
Let $\Gamma$ be a graph on $C$ such that $\{x,y\}\in E(\Gamma)$ if
and only if $x$ and $y$ are of the same color. Then $\phi$ is
precisely a ``constrained labeling'' in the sense of
Section~\ref{sec:coloring} and, in full agreement with
(\ref{eqn:chess-join}), $(\Delta_{r,2r-1})^{\ast (k+1)}$ is the
simplicial complex of all admissible labelings. In other words
$(\Delta_{r,2r-1})^{\ast (k+1)}$ is a well-chosen ``configuration
space'' \cite[Section 14.1]{Z04} associated to the colored
Tverberg problem. The associated ``test space'' is $W_r^{\oplus
d}$ and within the framework of ``Configuration Space -- Test
Map''-scheme \cite[Section 14.1]{Z04}, the problem is reduced to
the nonexistence of a $\mathbb{Z}/r$-equivariant map
(\ref{eqn:nasa}), see the original papers \cite{ZV92, VZ94} or
reviews  \cite{M,Z04, User2} for more detailed presentation.

\subsection{The breakthrough of Blagojevi\' c and Ziegler}

It is amusing to see how ingenious and astonishingly simple was
the new idea of Blagojevi\' c and Ziegler \cite{B-Z} leading to
the bound $t(r-1,d,d)=r-1$ for a prime $r$. They observed that if
$C\subset \mathbb{R}^d$ is a set of size $(r-1)(d+1)$ which is
evenly colored by $d+1$ colors, then it is natural to add one more
point $x\in \mathbb{R}^d$ and one more color (which corresponds to
the added point $x$). The enlarged set $C^+ = C\cup\{x\}$ is
colored by $d+2$ colors and it is natural to ask whether one can
find $r$ vertex-disjoint rainbow simplices which have a nonempty
intersection. Here, as in Section~\ref{sec:coloring}, a simplex is
rainbow if all its vertices are colored by different colors.

By using exactly the same translation as above, and in perfect
analogy with (\ref{eqn:nasa}), one is immediately led to the
question of the existence of a $\mathbb{Z}/r$-equivariant map
\begin{equation}
F : (\Delta_{r,r-1})^{\ast d}\ast [r] \rightarrow S(W_r^{\oplus
d}).
\end{equation}

\subsection{Examples}

Here are some examples of old and new colored Tverberg theorems,
rephrased as statements about simplicial maps of complexes
(graphs) $K_{p_1,p_2,\ldots, p_k}= [p_1]\ast [p_2]\ast\ldots\ast
[p_k]$. The connection with Problem~\ref{prob:color}
(Section~\ref{sec:CTP}) is established by an observation that a
coloring $C = C_1\cup\ldots\cup C_k$ of a set $C\subset
\mathbb{R}^d$, where $\vert C_i\vert = p_i$, defines a simplicial
map $\phi : K_{p_1,p_2,\ldots, p_k} \rightarrow \mathbb{R}^d$.
\begin{equation}
\label{K33} (K_{3,3} \longrightarrow {\mathbb R}^2) \Rightarrow
(2- {\rm crossing})
\end{equation}

\begin{equation}
\label{K333} (K_{3,3,3} {\longrightarrow} {\mathbb R}^2)
\Rightarrow (3- {\rm intersection})
\end{equation}

\begin{equation}
\label{K555} (K_{5,5,5} \longrightarrow {\mathbb R}^3) \Rightarrow
(3- {\rm crossing})
\end{equation}

\begin{equation}
\label{K4444} (K_{4,4,4,4} \longrightarrow {\mathbb R}^3)
\Rightarrow (4- {\rm intersection})
\end{equation}
The first of these results, claiming that for each (simplicial)
map $\phi : K_{3,3}\rightarrow \mathbb{R}^2$ there always exist
two intersecting vertex-disjoint edges in the image, is a
consequence of the nonplanarity of the complete bipartite graph
$K_{3,3}$. The second is an instance of a result of B\'{a}r\'{a}ny
and Larman \cite{BL}. It says that each collection of nine points
in the plane, evenly colored by three colors, can be partitioned
into three rainbow triangles which have a common point. A similar
conclusion have statement (\ref{K555}) which was in \cite{VZ94}
informally formulated as a statement about a constellation $5$
red, $5$ blue, and $5$ green stars in the outer space. Finally
(\ref{K4444}) is an instance of the result of Blagojevi\'{c},
Matschke, and Ziegler \cite[Corollary~2.4]{BMZ} saying that $4$
intersecting rainbow tetrahedra in $\mathbb{R}^3$ will always
appear if we are given sixteen points, evenly colored by four
colors.

All results (\ref{K33})--(\ref{K4444}) are best possible in the
sense that they provide exact values for the function $t(r,k,d)$.
All these results, possibly with exception of (\ref{K333}), remain
valid if $\phi$ is an arbitrary continuous map. Statements
(\ref{K333}) and (\ref{K4444}) are examples of the type A,  while
(\ref{K33}) and (\ref{K555}) are instances of type B colored
Tverberg theorem.

\section{Degrees of equivariant maps}\label{sec:main}

In this section we formulate our main results about equivariant
maps from joins of chessboard complexes. Short and elementary
proofs are given in Section~\ref{sec:proofs}.

Recall \cite[Section 2]{BLVZ} that $\Delta_{r,r-1}$ is an
orientable pseudomanifold. An associated fundamental homology
class is well-defined and a map $\phi : \Delta_{r,r-1}\rightarrow
M$ has a well-defined degree ${\rm deg}(\phi)$ for each orientable
$(r-2)$-dimensional manifold $M$. As before, $W_r$ is the
standard, $(r-1)$-dimensional real permutation representation of
$\mathbb{Z}/r$.

\begin{prop}\label{prop:degree}
The degree  ${\rm deg}(f)$ of each $\mathbb{Z}/r$-equivariant map
 \begin{equation}\label{eqn:degree}
f : (\Delta_{r,r-1})^{\ast d}\rightarrow S(W_r^{\oplus d})
 \end{equation}
is nonzero, provided $r$ is a prime number.  More precisely ${\rm
deg}(f)\equiv_{{\rm mod}\, r} (-1)^d$ and for each integer $m$
such that $m\equiv_{{\rm mod}\, r} (-1)^d$ there exists a
$\mathbb{Z}/r$-equivariant map $g : (\Delta_{r,r-1})^{\ast
d}\rightarrow S(W_r^{\oplus d})$ such that ${\rm deg}(g)=m$.
\end{prop}

Proposition~\ref{prop:degree} implies the following proposition
which establishes the main colored Tverberg type result of
\cite{B-Z}.

\begin{prop}\label{prop:posledica}
Let $r\geq 2$ be a prime and $d\geq 1$. Then there does not exist
a $\mathbb{Z}/r$-equivariant map
\begin{equation}\label{eqn:posledica}
F : (\Delta_{r,r-1})^{\ast d}\ast [r] \rightarrow S(W_r^{\oplus
d}).
\end{equation}
\end{prop}
The formal similarity of statements (\ref{eqn:nasa}),
(\ref{eqn:degree}), and (\ref{eqn:posledica}), and our general
emphasis on the equivariant maps from joins of chessboard
complexes, serve as a motivation for the following general result.

\begin{theo}\label{thm:general}
Suppose that $X$ is a $(\nu-1)$-connected, free
$\mathbb{Z}/r$-complex where $r$ is a prime number. Suppose that
$$U \cong W_r^{\oplus l}\oplus V $$
where $W_r$ is the standard $(r-1)$-dimensional permutation
representation of $\mathbb{Z}/r$ and $V$ an arbitrary real
fixed-point-free representation of dimension $\leq\nu$. Then there
does not exist a $\mathbb{Z}/r$-equivariant map
\begin{equation}\label{eqn:general}
f : (\Delta_{r,r-1})^{\ast l}\ast X \rightarrow S(W_r^{\oplus
l}\oplus V).
\end{equation}
\end{theo}
Using the known fact \cite{BLVZ} that $\Delta_{s,t}$ is
$(\nu-1)$-connected where $$\nu = {\rm
min}\{s,t,\lfloor\frac{1}{3}(s+t+1)\rfloor\}-1,$$
Theorem~\ref{thm:general} specializes to results claiming
nonexistence of equivariant maps of the form
\begin{equation}\label{eqn:gen-posledica}
f : (\Delta_{r,r-1})^{\ast l}\ast
\Delta_{r,s_1}\ldots\ast\Delta_{r,s_k} \rightarrow S(W_r^{\oplus
(d+1)})
\end{equation}
for an appropriate choice of parameters $s_1,\ldots,s_k$ and $l$,
which are carefully chosen to allow an application of
Theorem~\ref{thm:general}. Since $[r]$ is nothing but
$\Delta_{r,1}$ we observe that Proposition~\ref{prop:posledica} is
the simplest instance of (\ref{eqn:gen-posledica}).

The case $s_1=\ldots =s_l$ is of special interest and all this
together indicates that there should exist a plethora of colored
Tverberg results of mixed type A and type B in the sense of
\cite{Z04}.

\section{Colored Tverberg results of mixed type}

Here we specialize further and list the first consequences of
Theorem~\ref{thm:general}. We initially focus our attention to the
case $s_1=\ldots=s_k$.

\subsection{The case $s_1=\ldots=s_k=2r-1$}

Since $\Delta_{r,2r-1}$ is $(r-2)$-connected, the complex
$(\Delta_{r,2r-1})^{\ast k}$ is $(rk-2)$-connected. It follows
from Theorem~\ref{thm:general} that there does not exist a
$\mathbb{Z}/r$-equivariant map
\begin{equation}\label{eqn:gen-slucaj-1}
f : (\Delta_{r,r-1})^{\ast l}\ast (\Delta_{r,2r-1})^{\ast k}
\rightarrow S(W_r^{\oplus (d+1)})\cong S(W_r^{\oplus l})\ast
S(W_r^{\oplus (d-l+1)})
\end{equation}
provided
 \begin{equation}\label{eqn:nejednakost-1}
  (r-1)(d-l+1)+1\leq rk.
 \end{equation}
From here one immediately deduces the following proposition.
\begin{prop}\label{prop:slucaj-1}
Suppose that $C\subset \mathbb{R}^d$ is a collection of $N=(r-1)l
+ (2r-1)k$ points in $\mathbb{R}^d$ colored by $k+l$ colors and
let $C = \cup_{j=1}^{k+l}~C_j$ be the associated partition of $C$
into monochromatic parts. Assume that $\vert C_i\vert = r-1$ for
$i=1,\ldots,l$ and $\vert C_i\vert = 2r-1$ for $i=l+1,\ldots,l+k$.
Assume that the inequality {\em(\ref{eqn:nejednakost-1})}\! is
satisfied and that $r$ is a prime number. Then there exist $r$
vertex-disjoint rainbow simplices which have a nonempty
intersection.
\end{prop}
In the case $k=1$ we obtain the following result.

\begin{theo}\label{thm:k=1}
Suppose that $C_1,...,C_d, C_{d+1}$ are (monochromatic) sets in
$\mathbb{R}^d$ colored by $d+1$ distinct colors such that
$C_{d+1}$ has $2r-1$ elements while each of the remaining sets has
cardinality $r-1$. Then one can find $r$ vertex-disjoint rainbow
simplices with a nonempty intersection.
\end{theo}
Let us compare Theorem~\ref{thm:k=1} to the original result of
Blagojevi\'{c} and Ziegler \cite{B-Z}. Suppose one is interested
in conditions which guarantee the existence of $r$ intersecting
rainbow simplices. In Blagojevi\'{c}-Ziegler approach r+1 has to
be a prime and in our approach to the problem r is a prime number.
Neglecting for a moment this difference we observe that \cite{B-Z}
requires $r$ points in each of $d+1$ color classes whereas we ask
for $r-1$ points in $d$ color classes  and $2r-1$ points in the
remaining color class. The difference between the total numbers of
points in these two cases is $(r-1)d+2r-1-r(d+1)=r-(d+1)$. So, in
some (not direct) sense, their result gives more when $r$ is
greater than $d+1$, and our in the other case.

In the only case when both $r$ and $r+1$ are primes (the case
$r=2$) our approach yields a colorful extension of Radon's
theorem.
\begin{cor}\label{cor:Radon} (Colored Radon theorem) Let $C$ be a
collection of $d+3$ points in $\mathbb{R}^d$, three of the same
color and the remaining points all of different colors. Then there
exist two vertex-disjoint rainbow simplices which have a nonempty
intersection.
\end{cor}

Recall that the classical Radon theorem says that for each set $C$
of $d+2$ points in $\mathbb{R}^d$ there exist disjoint subsets
$C_1$ and $C_2$ of $C$ such that ${\rm conv}(C_1)\cap {\rm
conv}(C_2)\neq\emptyset$. Corollary~\ref{cor:Radon} says that if
we add one more point to $C$ and prescribe in advance a
three-element subset $D\subset C$, then there exist disjoint
subsets $C_1$ and $C_2$ of $C$ with intersecting convex hulls such
that the intersection $C_i\cap D$ is either empty or a singleton.

\medskip

The following corollary shows that Proposition~\ref{prop:slucaj-1}
is in some sense a mixed type A and type B colored Tverberg
theorem, \cite{Z04}.
\begin{cor}
If $l=0$ then the inequality {\em (\ref{eqn:nejednakost-1})}\!
reduces to $r\leq d/(d-k+1)$ and Proposition~\ref{prop:slucaj-1}
reduces to the type B colored Tverberg theorem, {\em \cite{VZ94,
Z04}}.
\end{cor}

\subsection{The case $r=2p-1$ and $s_1=\ldots=s_k=p$}

Suppose that $r=2p-1$ is an odd prime. It follows from
Theorem~\ref{thm:general} that there does not exist a
$\mathbb{Z}/r$-equivariant map
\begin{equation}\label{eqn:gen-slucaj-2}
f : (\Delta_{2p-1,2p-2})^{\ast l}\ast (\Delta_{2p-1,p})^{\ast k}
\rightarrow S(W_r^{\oplus (d+1)})\cong S(W_r^{\oplus l})\ast
S(W_r^{\oplus (d-l+1)})
\end{equation}
provided
 \begin{equation}\label{eqn:nejednakost-2}
  (r-1)(d-l+1)+1\leq pk.
 \end{equation}
\begin{prop}\label{prop:slucaj-2}
Suppose that $C\subset \mathbb{R}^d$ is a collection of $N=(r-1)l
+ pk$ points in $\mathbb{R}^d$ colored by $k+l$ colors and let $C
= \cup_{j=1}^{k+l}~C_j$ be the associated partition of $C$ into
monochromatic parts. Assume that $\vert C_i\vert = r-1$ for
$i=1,\ldots,l$ and $\vert C_i\vert = p$ for $i=l+1,\ldots,l+k$.
Assume that the inequality {\em (\ref{eqn:nejednakost-2})} is
satisfied and that $r$ is a prime number. Then there exist $r$
vertex-disjoint rainbow simplices which have a nonempty
intersection.
\end{prop}

Proposition~\ref{prop:slucaj-2} specializes for particular values
of parameters $r=2p-1,k,l,d$ to results that also deserve closer
inspection.

\medskip
For example if $l=0$ then the condition (\ref{eqn:nejednakost-2})
is fulfilled if we assume the equality $pk=(r-1)(d+1)+1$. Since
$(r-1)(d+1)+1$ is precisely the Tverberg number for
$r$-intersections in $d$-dimensional space, we observe that
Proposition~\ref{prop:slucaj-2} is also a refinement of the
classical (monochromatic) Tverberg theorem.

\begin{exam}
Choose $d=4,p=3,r=5,k=7$. Then Proposition~\ref{prop:slucaj-2}
says that if $21$ point in $\mathbb{R}^4$ is colored by $7$ colors
then there always exist $5$ vertex-disjoint rainbow simplices with
a nonempty intersection.
\end{exam}

\section{Proofs}\label{sec:proofs}
\subsection{Mapping degrees of equivariant maps}
The proof of Proposition~\ref{prop:degree} relies on a general
result (Proposition~\ref{prop:congruence}) which is an instance of
the typical ``comparison principle'' in the degree theory for
equivariant maps, \cite[page 4]{K-B}. This result can be also seen
as a relative of theorems about mapping degrees of equivariant
maps between representation spheres, see [Chap.\ II,
Proposition~4.12]\cite{tD} for an example, and \cite[page 139]{tD}
for a brief guide to other results of similar nature.

\begin{prop}\label{prop:congruence}
Suppose that $M$ is a triangulated, compact, orientable,
$n$-di\-men\-sio\-nal pseudomanifold. Let $G$ be a finite group
which acts freely and simplicially on $M$ and let $S(W)$ be a
$G$-invariant sphere in a real, $(n+1)$-dimensional
$G$-representation $W$. Suppose that $M$ and $S(W)$ have the same
orientation character, i.e.\ each element of $G$ either preserves
orientations of both $M$ and $S(W)$, or it reverses both of them.
Then for any two $G$-equivariant maps $f, g : M\rightarrow S(W)$,
\begin{equation}\label{eqn:kongruencija}
{\rm deg}(f) \equiv {\rm deg}(g) \quad {\rm mod}\, \vert G\vert.
\end{equation}
\end{prop}

\medskip\noindent
{\bf Proof:} Let $F : M\times I \rightarrow W$ be a
$G$-equivariant homotopy between maps $i \circ f$ and $i\circ g$
transverse to $0\in W$, where $i : S(W)\rightarrow W$ is the
inclusion map. Since the subspace $\Sigma\subset M$ of singular
points has dimension $\leq n-2$, the set $\Sigma\times I$ has
dimension $\leq n-1$, hence we can assume that $0\notin
F(\Sigma\times I)$.

It follows that the set $Z(F):=F^{-1}(0)$ is finite and consists
of nonsingular points. The set $Z(F)$ is clearly $G$-invariant.
For each $x\in Z(F)$ choose an open ball $V_x\ni x$ such that $V:=
\cup_{x\in Z(F)}~V_x$ is $G$-invariant and $V_x\cap V_y =
\emptyset$ for $x\neq y$. Let $S_x^n:=\partial(V_x)\cong S^n$ be
the boundary of $V_x$.

Let $N := (M\times I)\setminus V$, $M_0:= M\times \{0\}$ and
$M_1:= M\times \{1\}$. By construction there is a relation among
(properly oriented) fundamental classes,
\begin{equation}\label{eqn:homoloska-relacija}
  [M_1] - [M_0] = \sum_{x\in Z(F)} [S_x^n]
\end{equation}
inside the homology group $H_n(N,\mathbb{Z})$. The map $F_\ast :
H_n(N,\mathbb{Z})\rightarrow H_n(S(W),\mathbb{Z})$ maps the
relation (\ref{eqn:homoloska-relacija}) into the desired
congruence (\ref{eqn:kongruencija}). \hfill $\square$

\begin{rem}\label{rem:character}{\em
The condition in Proposition~\ref{prop:congruence} that $M$ and
$S(W)$ have the same orientation character is trivially fulfilled
if $G$ is a group with odd number of elements, in particular if
$G=\mathbb{Z}/r$ where $r$ is an odd prime, since a group of odd
order does not admit a nontrivial, one-dimensional real
representation.}
\end{rem}

\subsection{Canonical equivariant maps}\label{sec:canonical}

Proposition~\ref{prop:congruence} reduces the problem of
evaluating the $({\rm mod}\, r)$-degree of an arbitrary
$\mathbb{Z}/r$-equivariant map $f : M\rightarrow S(V)$ to the much
easier problem of testing a well chosen (canonical) map of this
kind.

\begin{defin} Let $\sigma^{m-1}$ be the simplex spanned by $[m]$ and
let $[m]^{(k)}:= \{A\subset [m]\mid \vert A\vert\leq k\}$ be its
$(k-1)$-skeleton, in particular
$[m]^{(m-1)}=\partial(\sigma^{m-1})$ is a triangulation of a
sphere. Define
 \begin{equation}\label{eqn:preslikavanje}
\xi=\xi_{m,k}: \Delta_{m,k}\rightarrow [m]^{(k)}
 \end{equation}
as the projection which sends a non-taking rook placement
$S=\{(i_1,j_1),\ldots, (i_p,j_p)\}\subset [m]\times [k]$ to the
set $\xi(S)=\{i_1,\ldots, i_p\}\subset [m]$.
\end{defin}

\begin{prop}\label{prop:mali-primer}
The degree of the map $\xi_{r,r-1} : \Delta_{r,r-1} \rightarrow
[r]^{(r-1)}$ is $${\rm deg}(\xi_{r,r-1})=(-1)^{r+1}(r-1)!$$
\end{prop}

\medskip\noindent
{\bf Proof:} Each simplex $\sigma\in\Delta_{r,r-1}$ can be
associated a unique permutation $\pi\in S_r$ such that
  \begin{equation}\label{eqn:fundamental}
 \sigma = \sigma_\pi = \{(\pi_1,1), (\pi_2,2),\ldots,
 (\pi_{r-1},r-1)\}.
 \end{equation}
It is not difficult to observe, cf.\ \cite[page 29]{BLVZ}, that if
${\widetilde\sigma}_\pi$ is the associated ordered simplex then a
fundamental class of $\Delta_{r,r-1}$ is represented by the
simplicial chain
$$
[\Delta_{r,r-1}] = \sum_{\pi\in S_r} (-1)^{{\rm sgn}(\pi)}
{\widetilde\sigma}_\pi.
$$
A fundamental class of $[r]^{(r-1)}=\partial(\sigma^{r-1})$ is
$$
[\partial(\sigma^{r-1})] = \sum_{i=1}^r (-1)^{i-1}
(1,\ldots,\widehat{i},\ldots,r).
$$
Since
$$\xi_\ast
(\widetilde{\sigma}_\pi)=(\pi_1,\pi_2,\ldots,\pi_{r-1},\widehat{\pi}_r)=(-1)^{{\rm
sgn}(\pi)+r-j}(1,\ldots,\widehat{j},\ldots,r)
$$
where $j:=\pi(r)$, we observe that
$$ \xi_\ast([\Delta_{r,r-1}]) =
(-1)^{r+1} (r-1)! [\partial(\sigma^{r-1})]
$$
which completes the proof of the proposition. \hfill $\square$

\begin{rem}{\rm
A more geometric proof of Proposition~\ref{prop:mali-primer} is
based on a simple algebraic count of points in the preimage
$\xi^{-1}(x)$ where $x$ is the barycenter of a top dimensional
simplex in $[r]^{(r-1)}$. Indeed, the ordered simplices mapped to
$(1,2,\ldots,r-1)$ are precisely the simplices of the form
$\tau_\alpha:=((1,\alpha_1),\ldots, (r-1,\alpha_{r-1}))$ for some
permutation $\alpha\in S_{r-1}$. It remains to be observed that
all simplices $\tau_\alpha$ have the same orientation.}
\end{rem}
The following corollary of the proof of
Proposition~\ref{prop:mali-primer} is not used
(Remark~\ref{rem:character}) in the proof of
Proposition~\ref{prop:degree} however it is a natural companion of
Proposition~\ref{prop:congruence}.
\begin{cor}
The pseudomanifolds $[\Delta_{r,r-1}]$ and
$[r]^{(r-1)}=\partial(\sigma^{r-1})$ have the same orientation
character with respect to the action of the symmetric group $S_r$.
\end{cor}

\subsection{Proofs}

\noindent {\bf Proof of Proposition~\ref{prop:degree}:} According
to \cite{BLVZ} the chessboard complex $\Delta_{r,r-1}$ is an
orientable pseudomanifold with a free action of the group
$\mathbb{Z}/r$. The same holds for the join
$(\Delta_{r,r-1})^{\ast d}$ so, in light of
Proposition~\ref{prop:congruence}, it is sufficient to exhibit a
canonical, $\mathbb{Z}/r$-equivariant map
$$
\pi : (\Delta_{r,r-1})^{\ast d}\rightarrow S(W_r^{\oplus d})
$$
with a known degree. Since $S(W_r)\cong [r]^{(r-1)}$ as
$\mathbb{Z}/r$-spaces, and
$$
S(W_r^{\oplus d})\cong S(W_r)^{\ast d}\cong ([r]^{(r-1)})^{\ast
d},
$$
we observe that $\xi := (\xi_{r,r-1})^{\ast d}$ is an example of
such a map. Hence the desired relation ${\rm
deg}(\xi)=[(r-1)!]^d\equiv_{{\rm mod}\, r} (-1)^d$ is a
consequence of  Proposition~\ref{prop:mali-primer}. \hfill
$\square$

\bigskip\noindent
{\bf Proof of Proposition~\ref{prop:posledica}:} Since the cone
$$
{\rm Cone}[(\Delta_{r,r-1})^{\ast d}]\cong (\Delta_{r,r-1})^{\ast
d}\ast [1]
$$
is a subcomplex of $(\Delta_{r,r-1})^{\ast d}\ast [r]$, we observe
that if an equivariant map
$$
F : (\Delta_{r,r-1})^{\ast d}\ast [r] \rightarrow S(W_r^{\oplus
d})
$$
exists, then its restriction on the subcomplex
$(\Delta_{r,r-1})^{\ast d}$ would have a zero degree, which is in
contradiction with Proposition~\ref{prop:degree}. \hfill $\square$

\medskip\noindent {\bf Proof of Theorem~\ref{thm:general}:}
Suppose that there exists a $\mathbb{Z}/r$-equivariant map
described in line (\ref{eqn:general}). Since $X$ is
$(\nu-1)$-connected and $S(V)$ is $(\nu-1)$-dimensional, there
exists a $\mathbb{Z}/r$-equivariant map $\alpha_1 :
S(V)\rightarrow X$. Consequently there exists a
$\mathbb{Z}/r$-equivariant map
\begin{equation}\label{eqn:nova-1}
\alpha = Id\ast\alpha_1  : (\Delta_{r,r-1})^{\ast d}\ast S(V)
\rightarrow (\Delta_{r,r-1})^{\ast d}\ast X
\end{equation}
and the composition map
\begin{equation}\label{eqn:nova-2}
f\circ\alpha : (\Delta_{r,r-1})^{\ast d}\ast S(V) \rightarrow
S(W_r^{\oplus d}\oplus V) = S(W_r^{\oplus d})\ast S(V).
\end{equation}
It follows from Proposition~\ref{prop:degree} and
Proposition~\ref{prop:congruence} that the degree of this map is
nonzero. On the other hand the degree of this map must be zero
since $S(V)$ is a $(\nu-1)$-dimensional sphere, $X$ is
$(\nu-1)$-connected and consequently the $(Id\circ\alpha_1)$-image
of the fundamental class of $(\Delta_{r,r-1})^{\ast d}\ast S(V)$
must be zero. This contradiction completes the proof of the
theorem. \hfill $\square$

\section{Combinatorial geometry on vector bundles}

{\em Combinatorial geometry on vector bundles} \cite{TV-bundle} is
a general program of extending combinatorial geometric results
about finite sets of points in $\mathbb{R}^d$, in particular the
theorems of Tverberg type, to the case of vector bundles, where
they become combinatorial geometric statements about finite
families of continuous cross-sections.

In the case of the canonical bundle over a Grassmann manifold,
these results include theorems about common affine $k$-dimensional
transversals of sets in $\mathbb{R}^d$. An example is the
following statement, a consequence of
\cite[Theorem~3.1]{TV-bundle},
$$
(K_{6,6} \rightarrow R^3) \Rightarrow (4 \mapsto {\rm line})
$$
which says that for every collection of $6$ red and $6$ blue
points in $\mathbb{R}^3$, there always exists a collection of four
vertex-disjoint edges with endpoints of different color (rainbow
edges) which a admit a common line transversal.

As demonstrated in \cite{TV-bundle}, the methods of {\em
parameterized ideal valued index theory} allow a systematic
approach to this problem, in particular all results of Tverberg
type formulated in earlier sections have their vector bundle
analogues. The reader is referred to \cite{BMZ-2} for very
interesting new results of this type, in the context of general
Tverberg-Vre\'{c}ica problem \cite{TV}.


\begin{thebibliography}{abcd}


\bibitem[Ata]{Ata04} C.A.~Athanasiadis. Decompositions and connectivity
of matching and chessboard complexes. \textit{Discrete Comput.\
Geom.} 31 (2004), 395--403.

\bibitem[A-F]{AuFie07}
S.~Ault, Z.~Fiedorowicz. Symmetric homology of algebras.
arXiv:0708.1575v54 [math.AT] 5 Nov 2007.

\bibitem[BFL]{BFL} I.~B\'{a}r\'{a}ny, Z.~F\"{u}redy, and L.~Lov\'{a}sz.
On the number of halving planes. \textit{Combinatorica},
10:175--183, 1990.

\bibitem[BL]{BL} I.~B\'{a}r\'{a}ny, D.G.~Larman. A colored version
of Tverberg's theorem. \textit{J.~London Math.\ Soc.}, II.\ Ser.,
45:314--320, 1992.

\bibitem[BSS]{BSS} I.~B\'{a}r\'{a}ny, S.B.~Shlosman, S.~Sz\" ucs.
On a topological generalization of a theorem of Tverberg.
\textit{J.~London Math.\ Soc.}, II.\ Ser., 23:158--164, 1981.


\bibitem[B-Z]{B-Z} P.V.M.~Blagojevi\' c, G.M.~Ziegler. Optimal
bounds for the colored Tverberg problem, arXiv:0910.4987v1
[math.CO].

\bibitem[BMZ]{BMZ} P.V.M.~Blagojevi\' c, B.~Matschke, G.M.~Ziegler.
Optimal bounds for the colored Tverberg problem, arXiv:0910.4987v2
[math.CO].

\bibitem[BMZ-2]{BMZ-2} P.V.M.~Blagojevi\' c, B.~Matschke, G.M.~Ziegler.
Optimal bounds for a colorful Tverberg--Vre\' cica problem,
arXiv:0911.2692v2 [math.AT].

\bibitem[BLVZ]{BLVZ} A.~Bj\" orner, L.~Lov\' asz, S.T.~Vre\' cica, and
R.T.\ \v Zivaljevi\' c. Chessboard complexes and matching
complexes. \textit{J.\ London Math.\ Soc.\ (2)} 49 (1994), 25--39.

\bibitem[tD]{tD} T.~tom Dieck. \textit{Transformation Groups}. de
Gruyter Studies in Mathematics 8, Berlin 1987.

\bibitem[Fie]{Fie07} Z.~Fiedorowicz. Question about a simplicial complex.
\textit{Algebraic Topology Discussion List} (maintained by Don
Davis), \url{http://www.lehigh.edu/~dmd1/zf93}.

\bibitem[F-H]{FrHa98} J.~Friedman and P.\ Hanlon. On the Betti numbers
of chessboard complexes. \textit{J.\ Algebraic Combin.} 8 (1998),
193--203.

\bibitem[G]{Ga79} P.F.~Garst. \textit{Cohen-Macaulay complexes and group
actions}. Ph.D.\ Thesis, Univ.\ of Wisconsin-Madison, 1979.


\bibitem[H]{H} S.~Hell. Tverberg's theorem with constraints.
\textit{J.\ Combinatorial Theory}, Ser.\ A 115:1402--1406, 2008.


\bibitem[J]{Jo-book} J.~Jonsson. \textit{Simplicial Complexes of
Graphs}. Lecture Notes in Mathematics, Vol.\ 1928. Springer 2008.

\bibitem[J1]{Jo07} J.~Jonsson. Exact sequences for the homology of the
matching complex. \textit{Journal of Combinatorial Theory}, Series
A 115 (2008), no.\ 8, 1504--1526.

\bibitem[J2]{Jo07-2} J.~Jonsson. On the $3$-torsion part of the
homology of the chessboard complex. \textit{Annals of
Combinatorics}, accepted.

\bibitem[K-B]{K-B} A.~Kushkuley and Z.I.~Balanov. \textit{Geometric Methods in
Degree Theory for Equi\-va\-riant Maps}. Lecture Notes in Math.\
1632, Springer, Berlin 1996.


\bibitem[M]{M} J.~Matou\v sek. \textit{Using the Borsuk-Ulam Theorem.
Lectures on Topological Methods in Combinatorics and Geometry}.
Universitext, Springer-Verlag, Heidelberg, 2003.

\bibitem[R-R]{ReRo00} V.~Reiner and J.~Roberts. Minimal resolutions and
homology of chessboard and matching complexes. \textit{J.\
Algebraic Combin.} 11 (2000), 135--154.

\bibitem[S-W]{ShaWa04} J.\ Shareshian and M.L.~Wachs. Torsion in the
matching complex and chessboard complex. \textit{Advances in
Mathematics} 212 (2007) 525--570.

\bibitem[T-V]{TV} H.\ Tverberg and S. Vre\'{c}ica. On generalizations of Radon's
theorem and the Ham sandwich theorem. \textit{European J.\ Comb.},
14:259–-264, 1993.

\bibitem[V\v Z94]{VZ94} S.~Vre\' cica and R.~\v Zivaljevi\' c. New cases of
the colored Tverberg theorem. In H.~Barcelo and G.~Kalai, editors,
\textit{Jerusalem Combinatorics '93}, pp. 325--334.

\bibitem[V\v Z09]{VZ07} S.~Vre\' cica and R.~\v Zivaljevi\' c.
Cycle-free chessboard complexes and symmetric homology of
algebras. \textit{European J.\ Combinatorics} 30 (2009) 542–-554.

\bibitem[X]{X} S.~Vre\' cica and R.~\v Zivaljevi\' c. Chessboard complexes
indomitable, preprint, November 12, 2009.

\bibitem[${\rm X}^+$]{X+} S.T.~Vre\' cica and R.T.~\v Zivaljevi\' c. Chessboard complexes
indomitable, arXiv:0911.3512v1 [math.CO], November 18, 2009.


\bibitem[W]{Wa03} M.L.~Wachs. Topology of matching, chessboard, and
general bounded degree graph complexes, Dedicated to the memory of
Gian-Carlo Rota. \textit{Algebra Universalis} 49 (2003), 345--385.

\bibitem[Z]{Zie94} G.M.~Ziegler. Shellability of chessboard
complexes. \textit{Israel J.\ Math.} 87 (1994), 97--110.

\bibitem[\v ZV92]{ZV92} R.T.~\v Zivaljevi\' c and S.T.~Vre\' cica. The
colored Tverberg's problem and complexes of injective functions.
\textit{J.\ Combin.\ Theory Ser.\ A} 61 (1992), 309--318.

\bibitem[\v Z96]{User1}
R.~\v Zivaljevi\' c.
\newblock User's guide to equivariant methods
in combinatorics.  \textit{Pu\-bli\-cations de l'Institut
Mathematique} (Beograd), 59(73), 114--130, 1996.

\bibitem[\v Z98]{User2}
R.~\v Zivaljevi\' c.
\newblock User's guide to equivariant methods
in combinatorics II. \textit{Publi\-cations de l'Institut
Mathematique} (Beograd), 64(78) 1998, 107--132.

\bibitem[\v Z99]{TV-bundle} R.~\v Zivaljevi\'{c}.
The Tverberg-Vre\'{c}ica problem and the combinatorial geometry on
vector bundles. \textit{Israel J.\ Math}, 111:53–-76, 1999.


\bibitem[\v Z04]{Z04}
R.T. \v Zivaljevi\'{c}. Topological methods. Chapter 14 in
\textit{Handbook of Discrete and Computational Geometry}, J.E.\
Goodman, J.\ O'Rourke, eds, Chapman \& Hall/CRC 2004, 305--330.

\end{thebibliography}
\end{document}